\documentclass[final,1p,times]{elsarticle}
\usepackage{amsmath}
\usepackage{amssymb}
\usepackage{url}
\usepackage{amsfonts}
\usepackage{graphicx}
\usepackage{rotating}
\usepackage{floatflt,epsfig}
\usepackage{lineno,hyperref}
\usepackage{enumerate}
\usepackage{colortbl}
\usepackage{array,tabularx,tabulary,booktabs}
\usepackage{longtable}
\usepackage{multirow}
\usepackage{wrapfig}
\usepackage{subcaption}
\newcolumntype{^}{>{\currentrowstyle}}

\makeatletter
\def\ps@pprintTitle{%
   \let\@oddhead\@empty
   \let\@evenhead\@empty
   \let\@oddfoot\@empty
   \let\@evenfoot\@oddfoot
}
\makeatother

\newtheorem{lemma}{Lemma}

\newtheorem{proposition}{Proposition}

\newtheorem{construction}{Construction}

\newcommand{\proof}{\medskip\noindent{\bf Proof.~}}
\begin{document}
\renewcommand{\abstractname}{Abstract}
\renewcommand{\refname}{References}
\renewcommand{\tablename}{Table}
\renewcommand{\arraystretch}{0.9}
\thispagestyle{empty}
\sloppy

\begin{frontmatter}
\title{Classification of divisible design graphs with at most 39 vertices}

\author[01,02]{Dmitry~Panasenko}
\ead{makare95@mail.ru}

\author[01]{Leonid~Shalaginov}
\ead{44sh@mail.ru}

\address[01]{Chelyabinsk State University, Brat'ev Kashirinyh st. 129, Chelyabinsk, 454021, Russia}

\address[02]{Krasovskii Institute of Mathematics and Mechanics, S. Kovalevskaja st. 16, Yekaterinburg, 620990, Russia}

\begin{abstract}
A $k$-regular graph is called a divisible design graph (DDG for short) if its vertex set can be partitioned into $m$ classes of size $n$, such that two distinct vertices from the same class have exactly $\lambda_1$ common neighbours, and two vertices from different classes have exactly $\lambda_2$ common neighbours. A DDG with $m = 1$, $n = 1$, or $\lambda_1 = \lambda_2$ is called improper, otherwise it is called proper. We present new constructions of DDGs and, using a computer enumeration algorithm, we find all proper connected DDGs with at most 39 vertices, except for three tuples of parameters: $(32,15,6,7,4,8)$, $(32,17,8,9,4,8)$, $(36,24,15,16,4,9)$.
\end{abstract}

\begin{keyword} divisible design graph \sep divisible design \sep walk-regular graph
\vspace{\baselineskip}
\MSC[2010] 05C50\sep 05E10\sep 15A18
\end{keyword}
\end{frontmatter}

\section{Introduction}

An incidence structure with constant block size $k$ is a \emph{(group) divisible design} whenever the point set can be partitioned into $m$ classes of size $n$, such that two points from one class occur together in exactly $\lambda_1$ blocks, and two points from different classes occur together in exactly $\lambda_2$ blocks. A divisible design $D$ is called \emph{symmetric} or to have the \emph{dual property} (SDD for short) if the dual of $D$ (that is, the design with the transposed incidence matrix) is again a divisible design with the same parameters as $D$. A \emph{divisible design graph} is a graph whose adjacency matrix is the incidence matrix of a symmetric divisible design. A DDG with $m = 1$, $n = 1$, or $\lambda_1 = \lambda_2$ is called improper (these DDGs are $(v,k,\lambda)$-graphs), otherwise it is called proper. 

At first, divisible design graphs were studied in master's thesis by M.A.~Meulenberg \cite{M2008} and then studied in more detail in 2011 in the following paper by W.H.~Haemers, H.~Kharaghani and M.A.~Meulenberg \cite{HKM2011} and in 2011-2013 in two following papers by D.~Crnkovic and W.H.~Haemers \cite{CH2011, CH2014}.

In 2008 M.A.~Meulenberg presented the list of feasible parameters of proper DDGs up to 50 vertices. In 2011-2013 feasible parameters of proper DDGs up to 27 vertices were studied and the existence of graphs was resolved in all but one case, however, the exact number of graphs corresponding to these tuples of parameters remained unknown.

In this paper we present new constructions of DDGs and find all proper connected DDGs with at most 39 vertices, except for three tuples of parameters: $(32,15,6,7,4,8)$, $(32,17,8,9,4,8)$, $(36,24,15,16,4,9)$. 

The paper is organised as follows. In Section 2 we give some definitions, notations and preliminaries about DDGs. In Section 3 we give an overview of known constructions of DDGs and in Section 4 we describe some new constructions. In Section 5 we describe some new sporadic constructions of DDGs. In Section 6 we describe the algorithm used for enumerating DDGs and in Section 7 we present the results.

\section{Preliminaries}

Let $I_t$ and $J_t$ be an identity $t \times t$ and all-ones $t \times t$ matrix, respectively, and $K_{(m,n)}$ be $I_m \otimes J_n = diag(J_n, \ldots, J_n)$. Then a graph $\Gamma$ is a DDG with parameters $(v, k, \lambda_1, \lambda_2, m, n)$ if and only if $\Gamma$ has an adjacency matrix $A$ that satisfies

$$A^2 = kI_v + \lambda_1(K_{(m,n)} - I_v) + \lambda_2(J_v - K_{(m,n)}).$$

In DDGs $v = mn$, and taking row sums on both sides of the equation above gives

$$k^2 = k + \lambda_1(n - 1) + \lambda_2n(m - 1).$$

The formula for $A^2$ also gives us strong information about the eigenvalues of $A$ and their multiplicities.

\begin{lemma}[{\cite[Lemma 2.1]{HKM2011}}]\label{l1}
$A$ has at most five distinct eigenvalues $k$, $\sqrt{k - \lambda_1}$, $-\sqrt{k - \lambda_1}$, $\sqrt{k^2 - \lambda_2v}$, $-\sqrt{k^2 - \lambda_2v}$ with corresponding multiplicities $1, f_1, f_2, g_1, g_2$, where $f_1 + f_2 = m(n - 1)$ and $g_1 + g_2 = m - 1$. 
\end{lemma}

Some of the multiplicities may be 0 and some values may coincide. In general, the multiplicities $f_1, f_2, g_1$ and $g_2$ are not determined by the parameters, but if we know one, we can compute the rest because $f_1 + f_2 = m(n - 1), g_1 + g_2 = m - 1$ and $$\text{trace}(A) = 0 = k + (f_1 - f_2)\sqrt{k - \lambda_1} + (g_1 - g_2)\sqrt{k^2 - \lambda_2v}.$$

This equation leads to the following result.

\begin{lemma}[{\cite[Theorem 2.2]{HKM2011}}]\label{l2}
Consider a proper DDG with parameters $(v, k, \lambda_1, \lambda_2, m, n)$ and eigenvalue multiplicities $f_1, f_2, g_1, g_2$. Then

(1) $k - \lambda_1$ or $k^2 - \lambda_2v$ is a nonzero square.

(2) If $k - \lambda_1$ is not a square, then $f_1 = f_2 =m(n -1)/2$.

(3) If $k^2 - \lambda_2v$ is not a square, then $g_1 = g_2 = (m -1)/2$.
\end{lemma}

Let $V_1 \cup  V_2 \cup \ldots \cup V_t$ be the partition of the vertex set of a graph $\Gamma$ with the property that every vertex of $V_i$ has exactly $r_{ij}$ neighbours in $V_j$. Then $V_1 \cup  V_2 \cup \ldots \cup V_t$ will be an \emph{equitable $t$-partition} of $\Gamma$. Matrix $R = (r_{ij})_{t \times t}$ is called the \emph{quotient matrix} of the equitable partition.

\begin{lemma}[{\cite[Theorem 3.1]{HKM2011}}]\label{l3}
The vertex partition from the definition of a DDG (the canonical partition) is equitable and the quotient matrix $R$ satisfies $$R^2 = RR^T = (k^2 - \lambda_2v)I_m +\lambda_2nJ_m.$$

The eigenvalues of $R$ are $k$, $\sqrt{k^2 - \lambda_2v}$, $-\sqrt{k^2 - \lambda_2v}$ with corresponding multiplicities $1$, $g_1$, $g_2$.
\end{lemma}

\begin{lemma}[{\cite[Proposition 3.2]{HKM2011}}]\label{l4}
The quotient matrix $R$ of a DDG satisfies

$\sum\limits_{i} r_{ij} = k \text{ for } j = 1, \ldots, m,$

$\sum\limits_{i,j} r_{ij}^2 = \text{trace}(R^2) = m k^2 - (m-1)\lambda_2v,$

$0 \leqslant \text{trace}(R) = k + (g_1 - g_2)\sqrt{k^2 - \lambda_2v} \leqslant m(n-1).$
\end{lemma}

A graph is \emph{walk-regular}, whenever for every $l \geqslant 2$ the number of closed walks of length $l$ at a vertex $x$ is independent of the choice of $x$. Note that walk-regularity implies regularity (take $l = 2$). 

A DDG with four distinct eigenvalues is walk-regular, provided it is connected (\cite[Corollary 4.2]{CH2011}). A DDG with five distinct eigenvalues can also be walk-regular. To decide on this the following lemma can be used.

\begin{lemma}[{\cite[Theorem 4.3]{CH2011}}]\label{l5}
A proper DDG is walk-regular if and only if the quotient matrix $R$ has constant diagonal.
\end{lemma}

\section{Known constructions}

\subsection{$(v,k,\lambda)$-graphs and designs}

The \emph{incidence graph} of a design with incidence matrix $N$ is a bipartite graph with adjacency matrix

~

$\begin{bmatrix}
O & N\\
N^\top & O
\end{bmatrix}$.

\begin{construction}[{\cite[Construction 4.1]{HKM2011}}]\label{c1}
The incidence graph of an $(n,k,\lambda_1)$-design with ~~~~~~~~ $1 < k \leqslant n$ is a proper DDG with $\lambda_2 = 0$.
\end{construction}

\begin{proposition}[{\cite[Proposition 4.3]{HKM2011}}]\label{p1}
For a proper connected DDG $\Gamma$ the following are equivalent.

(1) $\lambda_2 = 0$.

(2) $\Gamma$ comes from Construction 1.
\end{proposition}

\begin{construction}[{\cite[Construction 4.4]{HKM2011}}]\label{c2}
If $A'$ is the adjacency matrix of an $(m,k',\lambda')$-graph $(1 \leqslant k' < n)$, then $A = A' \otimes J_n$ is the
adjacency matrix of a proper DDG with $k = \lambda_1 = nk'$, $\lambda_1 = n\lambda'$.
\end{construction}

\begin{proposition}[{\cite[Proposition 4.5]{HKM2011}}]\label{p2}
For a proper DDG $\Gamma$ the following are equivalent.

(1) $\lambda_1 = k$.

(2) $\Gamma$ comes from Construction 2.
\end{proposition}

\begin{construction}[{\cite[Construction 4.6]{HKM2011}}]\label{c3}
Let $A_1, \ldots, A_m$ $(m \geqslant 2)$ be the adjacency matrices of $m$ $(n,k',\lambda')$-graphs with $0 \leqslant k \leqslant n-2$. Then $A = J_v - K_{(m,n)} + diag(A_1, \ldots, A_m)$ is the adjacency matrix of a proper DDG with $k = k' + n(m-1)$, $\lambda_1 = \lambda' + n(m-1)$, $\lambda_2 = 2k - v$.
\end{construction}

\begin{proposition}[{\cite[Proposition 4.7]{HKM2011}}]\label{p3}
For a proper DDG $\Gamma$ the following are equivalent.

(1) $\lambda_2 = 2k - v$.

(2) $\Gamma$ comes from Construction 3.
\end{proposition}

\subsection{DDGs with $\lambda_1 = k - 1$}

The \emph{lexicographic product} or \emph{graph composition} $G[H]$ of graphs $G$ and $H$ is a graph such that the vertex set of $G[H]$ is $V(G) \times V(H)$ and adjacency defined by 
$$(u_1, u_2) \sim (v_1, v_2) \text{ if and only if } u_1 \sim v_1 \text{ or } (u_1 = v_1 \text{ and } u_2 \sim v_2).$$

\begin{construction}[{\cite[Theorem 4.11]{HKM2011}}]\label{c4}
If $G$ is a strongly regular graph with parameters $(v,k,\lambda, \lambda + 1)$, then $G[K_2]$ is a DDG with parameters $(2v, 2k + 1, 2k, 2\lambda+2, v, 2)$. If $G$ is $K_{x \times y}$, the complete multipartite graph containing $x$ parts of $y$ vertices, then $G[K_2]$ is a DDG with parameters $(2xy, 2y(x - 1) + 1, 2y(x - 1), 2y(x - 2) + 2, x, 2y)$.
\end{construction}

An involutive automorphism of a graph is called \emph{Seidel automorphism} if it interchanges only non-adjacent vertices. Permuting the rows (and not the columns) of the adjacency matrix of a graph according to Seidel automorphism is called \emph{dual Seidel switching}.

\begin{construction}[{\cite[Construction 2]{GHKS2019}}]\label{c5}
Let $\Gamma$ be a DDG obtained with the first case of Construction 4. Let $M$ be the adjacency matrix of $\Gamma$, and $P$ be a non-identity permutation matrix of the same size. Then $PM$ is the adjacency matrix of a DDG if and only if $P$ represents a Seidel automorphism.
\end{construction}

\begin{proposition}[{\cite[Theorem 2]{GHKS2019}}]\label{p4}
For a proper DDG $\Gamma$ with $\lambda_2 \not= 0$ the following are equivalent.

(1) $\lambda_1 = k - 1$.

(2) $\Gamma$ comes from Construction 4 or 5.
\end{proposition}

\subsection{Distance-regular graphs}

Suppose that $G$ is a graph with diameter $d$. For any vertex $u$ and for any integer $i$, where $0 \leqslant i \leqslant d$, let $G_i(u)$ denote the set of vertices at distance $i$ from $u$. If $u' \in G_i(u)$ and $w$ is a neighbour of $u'$, then $w$ must be at distance $i-1$, $i$ or $i+1$ from $u$. Let $c_i$, $a_i$ and $b_i$ denote the number of such vertices $w$. $G$ is a \emph{distance-regular graph} if and only if the parameters $c_i$, $a_i$, $b_i$ depend only on the distance $i$, and not on the choice of $u$ and $u'$ (i.e. $a_i+b_i+c_i=k=b_0, c_1=1$). The array $\{k, b_1, \ldots, b_{d-1}; 1, c_2, \ldots, c_d\}$ is called the \emph{intersection array} of the distance-regular graph. A distance-regular graph of diameter $d$ is called \emph{antipodal} if being at distance $d$ or $0$ defines an equivalence relation on the vertices. For a distance-regular graph the parameters $\lambda$ and $\mu$ give the number of common neighbours of a pair of vertices at distance 1 and 2, respectively (i.e. $\lambda=a_1, \mu=c_2$).

Distance regular graphs of diameter 2 are strongly regular graphs with parameters $(v, k, \lambda, \mu)$.

\begin{construction}[{\cite[Theorem 4.13]{HKM2011}}]\label{c6}
Suppose $\Gamma$ is an antipodal distance-regular graph of diameter 3. If $\lambda = \mu$, then $\Gamma$ is a proper DDG with parameters $(n(\mu n + 2),\mu n + 1, 0, \mu, \mu n + 2,n)$. If $\lambda = \mu - 2$, then the complement of $\Gamma$ is a proper DDG with parameters $(\mu n^2,\mu n(n-1),$ $\mu n(n-2),\mu(n-1)^2,\mu n,n)$.
\end{construction}

\subsection{Partial complements}

The \emph{partial complement} of a DDG is a graph whose adjacency matrix can be obtained as the complement of all blocks of the canonical partition except the diagonal blocks.

\begin{proposition}[{\cite[Proposition 4.15]{HKM2011}}]\label{p5}
The partial complement of a proper DDG $\Gamma$ is again a DDG if one of the following holds:

(1) The quotient matrix $R$ equals $t(J_m - I_m)$ for some $t \in \{1, \ldots , n-1\}$.

(2) $m = 2$.
\end{proposition}

Let $G$ be a $k$-regular graph on $v$ vertices with the smallest eigenvalue $\lambda_{min}$. A \emph{Hoffman coloring} of $G$ is a partition of the vertices into \emph{Hoffman-cocliques}, that is, cocliques meeting the Hoffman upper bound $c = v\lambda_{min}/(\lambda_{min} - k)$. An equitable partition of a $(v,k,\lambda)$-graph that satisfies (1) from Proposition 5 is a Hoffman coloring. 

\begin{construction}[{\cite[Construction 4.16]{HKM2011}}]\label{c7}
Let $\Gamma$ be a $(v,k,\lambda)$-graph. If $\Gamma$ has a Hoffman coloring or an equitable partition into two parts of equal size, then the partial complement is a DDG.
\end{construction}

\subsection{Hadamard matrices}

An $m \times m$ matrix $H$ is a \emph{Hadamard matrix} if every entry is $1$ or $-1$ and $HH^\top$ = $mI_m$. A Hadamard matrix $H$ is called \emph{graphical} if $H$ is symmetric with constant diagonal, and \emph{regular} if all row and column sums are equal.

\begin{construction}[{\cite[Construction 4.8]{HKM2011}}]\label{c8} Consider a regular graphical Hadamard matrix $H$ of order $m \geqslant 4$ and row sum $l = \pm \sqrt{m}$. Let $n \geqslant 2$. Replace each entry with value $-1$ by $J_n - I_n$, and each $+1$ by $I_n$, then we obtain the adjacency matrix of a DDG with parameters $(mn, n(m - l)/2 + l, (n - 2)(m - l)/2, n(m - 2l)/4 + l, m, n)$.
\end{construction}

\begin{construction}[{\cite[Construction 4.9]{HKM2011}}]\label{c9}Consider a regular graphical Hadamard matrix $H$ of order $l^2 \geqslant 4$ with diagonal entries $-1$ and row sum $l$. The graph with adjacency matrix

$A = \begin{bmatrix}
M & N & O\\
N & O & M\\
O & M & N
\end{bmatrix},$

where

$M = \displaystyle \frac{1}{2} \begin{bmatrix}
J_{l^2} + H & J_{l^2} + H\\
J_{l^2} + H & J_{l^2} + H
\end{bmatrix} \text{ and }
N = \displaystyle \frac{1}{2} \begin{bmatrix}
J_{l^2} + H & J_{l^2} - H\\
J_{l^2} - H & J_{l^2} + H
\end{bmatrix},$

is a DDG with parameters ($6l^2$, $2l^2 + l$, $l^2 + l$, $(l^2 + l)/2$, $3$, $2l^2$).
\end{construction}

\begin{construction}[{\cite[Theorem 3.2]{CH2011}}]\label{c10}
If there exist a regular graphical Hadamard matrix of order $4u^2$ with row sum $2u$ and a Hadamard matrix of order $2u^2$, then there exists a DDG with parameters
$(24u^2, 12u^2 - 2u, 4u^2 - 2u, 6u^2 - 2u, 12u^2, 2)$.
\end{construction}

\subsection{DDGs with parameters $(4n,n+2,n-2,2,4,n)$ and $(4n,3n-2,3n-6,2n-2,4,n)$}

An {\it $m\times n$-lattice graph} is a line  graph of complete bipartite graph $K_{m,n}$. 

\begin{construction}[{\cite[Theorem 1]{S2021}}]\label{c11}
An $m\times n$-lattice graph is a DDG if and only if $n=4$. These graphs have parameters $(4n,n+2,n-2,2,4,n)$. If $m=4$, then the graph is strongly regular with parameters $(16,6,2,2)$.
\end{construction}

\begin{construction}[{\cite[{Construction 4}]{S2021}}]\label{c12}
Let $M$ be the adjacency matrix of a $4\times n$-lattice graph. Let $M= \begin{bmatrix} M_{11} & M_{12}\\ M_{21}& M_{22}\\ \end{bmatrix}$, 
such that $M_{11}$ is the adjacency matrix of the subgraph $H$, which is isomorphic to the $2\times n$-lattice graph. Consider the permutation matrix $P = \begin{bmatrix} P_{11} & 0\\ 0 & I\\ \end{bmatrix}$,
where $P_{11}$ is the permutation matrix of the Seidel automorphism $\varphi$ corresponding to the central symmetry of $H$. Then $\begin{bmatrix} P_{11}M_{11} & M_{12}\\ M_{21}& M_{22}\\ \end{bmatrix}$ is the adjacency matrix of a DDG with parameters $(4n,n+2,n-2,2,4,n)$.
\end{construction}

\emph{The switching of edges} between two sets of vertices of a graph is the reversion of the adjacency of each pair of vertices, one from the first set and other from the second set. Thus, the edge set is changed so that an adjacent pair becomes nonadjacent and a nonadjacent pair becomes adjacent.

\begin{construction}[{\cite[Construction 5]{S2021}}]\label{c13}
Consider $C_{4t}[\overline{K_2}]$ ($t\geqslant 1$) and the $4$-cube, where $C_{4t}$ is the $4t$-cycle. Each of these graphs has an equitable partition with quotient matrix $J_4$. Consider some copies of $C_{4t}[\overline{K_2}]$ and some copies of the $4$-cube with a fixed equitable partition. Denote by $V_1, V_2, V_3,V_4$ the classes of this partition. The switching of edges between $V_1$ and $V_2$ and also between $V_3$ and $V_4$ gives a DDG with parameters $(4n,n+2,n-2,2,4,n)$. 
\end{construction}

For more information about equitable partitions of the $4$-cube and $C_{4t}[\overline{K_2}]$ see Lemmas 5-8 from \cite{S2021}.

Note that graphs obtained from Construction \ref{c11}, Construction \ref{c12} and Construction \ref{c13} have the same parameters, but different spectra. 

\begin{proposition}[{\cite[Theorem 1]{S2021}}]\label{p6}
Let $\Gamma$ be a DDG with parameters $(4n,n+2,n-2,2,4,n)$. 
Then:

(1) If $n$ is odd, then $\Gamma$ is isomorphic to $4\times n$-lattice graph.

(2) If $n$ is even, then $\Gamma$ comes from Construction \ref{c11}, \ref{c12} or \ref{c13}.
\end{proposition}

\begin{construction}[{\cite[Theorem 2]{S2021}}]\label{c14}
Let $\Gamma$ be a DDG with parameters $(4n,n+2,n-2,2,4,n)$.
The switching of edges between the union of two classes of the canonical partition and the remaining vertices gives a DDG with parameters $(4n,3n-2,3n-6,2n-2,4,n)$. 
\end{construction}

\section{New constructions}

\begin{construction}\label{c15}
Suppose $\Gamma$ is an antipodal distance-regular graph of diameter 3 with antipodal classes of size $r$. Denote by $A_i$ the matrix of a relation `to be at distance i' on the vertices of $\Gamma$. If $\lambda = \mu + 2$, then $A=A_1+A_3$ is the adjacency matrix of a DDG with parameters ($r(2\mu + 4)$, $2\mu + r + 2$, $r - 2$, $\mu + 2$, $2\mu + 4$, $r$).
\end{construction}

\proof The intersection array of an antipodal distance-regular graph of diameter $3$ is $\{k$, $\mu(r-1)$, $1$; $1$, $\mu$, $k\}$ \cite[p. 431]{BCN1989}. From this the statement of the construction follows straightforwardly. \hfill $\square$

\medskip

\begin{construction}\label{c16}
Suppose $\Gamma$ is a strongly regular graph with parameters $(v, k, \mu + 2, \mu)$ and has a Hoffman coloring with Hoffman-cocliques of size $n$ ($v = mn$). Let $A$ be an adjacency matrix of $\Gamma$, in which Hoffman-cocliques are located on the main diagonal, $K = K_{(m, n)}$, $I = I_v$. Then $A + K - I$ is the adjacency matrix of a DDG with parameters ($mn$, $k + n - 1$, $n + \mu - 2$, $\mu + \frac{2k}{m-1}$, $m$, $n$).
\end{construction}

\proof Let's calculate $(A + K - I)^2$ (we denote $J = J_v$).

$$(A + K - I)^2 = A^2 + K^2 + I^2 + AK + KA - 2A - 2K.$$

Since $A^2 = kI + \lambda A + \mu(J - I - A)$, $K^2 = nK$ and $AK = KA = {\textstyle \frac{k}{m-1}(J - K)}$, 

$$\begin{aligned}
(A + K - I)^2 = kI + \lambda A + \mu(J - I - A) + nK + I + {\textstyle \frac{2k}{m-1}(J - K)} - 2A - 2K = \\ = (k - \mu + 1)I + (\mu + {\textstyle \frac{2k}{m-1}})J + (n - {\textstyle \frac{2k}{m-1}} - 2)K + (\lambda - \mu - 2)A \text{.~}
\end{aligned}$$

Since $\lambda = \mu + 2$, 

$$(A + K - I)^2 = (k - \mu + 1)I + (\mu + {\textstyle \frac{2k}{m-1}})J + (n - {\textstyle\frac{2k}{m-1}} - 2)K,$$
 
which can be rewritten as  

$$(A + K - I)^2 = (k + n - 1)I + (n + \mu - 2)(K-I) + (\mu + {\textstyle\frac{2k}{m-1}})(J - K).$$

So, $A + K - I$ is the adjacency matrix of a DDG with parameters ($mn$, $k + n - 1$, $n + \mu - 2$, $\mu + {\textstyle\frac{2k}{m-1}}$, $m$, $n$). \hfill $\square$

\medskip

\begin{construction}\label{c17}
Let $\Gamma$ be a DDG with parameters $(v, k, \lambda_1, \lambda_2, m, n)$ and quotient matrix ~~~~~~~~~~~ $R = a I_m + b (J_m - I_m)$. Take $s$ copies of $\Gamma$ and label all blocks of the canonical partition in each copy with numbers $1, \ldots, m$. Then connect all vertices from blocks with the same label (adjacency inside the block does not change). The resulting graph is a DDG with parameters $(vs, k + (s-1)n, \lambda_1 + (s-1)n, \lambda_2, ms, n)$, if $\lambda_2 = 2b = 2a + (s-2)n$.
\end{construction}

\proof Consider two vertices from the same block of the same copy of $\Gamma$. Then they have $\lambda_1$ common neighbours in $\Gamma$ plus $n$ common neighbours in each other copy of $\Gamma$. Thus, the number of common neighbours equals $\lambda_1 + (s-1)n$. 

Consider two vertices from different blocks of the same copy of $\Gamma$. Then they have $\lambda_2$ common neighbours in $\Gamma$ and no new common neighbours in other copies of $\Gamma$. Thus, the number of common neighbours equals $\lambda_2$.

Consider two vertices from the same block of different copies of $\Gamma$. Then they have $a$ common neighbours in their blocks ($2a$ in total) plus $n$ common neighbours in each other copy of $\Gamma$. Thus, the number of common neighbours equals $2a + (s-2)n$.

Consider two vertices from different blocks of different copies of $\Gamma$. Then they have $b$ common neighbours in their copies of $\Gamma$ ($2b$ in total) and no new common neighbours in other copies of $\Gamma$. Thus, the number of common neighbours equals $2b$.

So, the resulting graph is a DDG with parameters $(vs, k + (s-1)n, \lambda_1 + (s-1)n, \lambda_2, ms, n)$, if $\lambda_2 = 2b = 2a + (s-2)n$. \hfill $\square$

\begin{construction}\label{c18}
Let $D$ be a symmetric divisible design with parameters $(v, k, \lambda_1,\lambda_2, m, n)$, such that every block contains $\frac{k}{m}$ points from each class. $\Gamma$ is the incidence graphs of $D$. Construct a new graph $\Gamma^{\ast}$ with the same vertex set as $\Gamma$, where vertices $x, y$ are adjacent in $\Gamma^{\ast}$ when $x, y$ are adjacent in $\Gamma$, or $x, y$ are points from different classes of $D$, or $x, y$ are blocks from different classes of the dual design of $D$. Then $\Gamma^{\ast}$ is a DDG with parameters ($2v, k + (m-1)n$, $\lambda_1 + (m-1)n$, $\lambda_2 + (m-2)n$, $2m$, $n$), if $2(m-1)\frac{k}{m} = \lambda_2 + (m-2)n$.
\end{construction}

\proof Consider two vertices, corresponding to two points from the same class of $D$. They have $\lambda_1$ common neighbours in $G$ plus $n$ common neighbours in each other class of points. Thus, we have $\lambda_1 + (m-1)n$ common neighbours. The proof goes similarly for the case when we consider two vertices corresponding to two blocks from the same class. 

Consider two vertices, corresponding to two points from different classes of $D$. They have $\lambda_2$ common neighbours in $G$ plus $n$ common neighbours in each other class of points. Thus, we have $\lambda_2 + (m-2)n$ common neighbours. The proof goes similarly for the case when we consider two vertices corresponding to two blocks from different classes.

Consider two vertices, the first corresponding to a point and the second corresponding to a block of $D$. They have no common neighbours in $G$. They have $\frac{k}{m}$ common neighbours in each class of points except for the class the first vertex belongs to. They also have $\frac{k}{m}$ common neighbours in each class of blocks except for the class the second vertex belongs to. Thus, we have $2(m-1)\frac{k}{m}$ common neighbours. 

So, $\Gamma^{\ast}$ is a DDG with parameters ($2v, k + (m-1)n$, $\lambda_1 + (m-1)n$, $\lambda_2 + (m-2)n$, $2m$, $n$), if $2(m-1)\frac{k}{m} = \lambda_2 + (m-2)n$. \hfill $\square$

\begin{construction}\label{c19}
Let $D$ be a symmetric divisible design with parameters $(v, k, \lambda_1,\lambda_2, m, n)$ with even $m$, such that every block contains $\frac{k}{m}$ points from each class. Let $\Gamma$ be the incidence graphs of $D$. Split the set of classes of blocks and the set of classes of points into pairs. Construct a new graph $\Gamma^{\ast}$ with the same vertex set as $\Gamma$, where vertices $x, y$ are adjacent in $\Gamma^{\ast}$ when $x, y$ are adjacent in $\Gamma$ or $x,y$ are from different classes of the same pair of classes partition. Then $\Gamma^{\ast}$ is a DDG with parameters ($2v, k + n$, $\lambda_1 + n$, $\lambda_2$, $2m$, $n$), if $\frac{2k}{m} = \lambda_2$.
\end{construction}

\proof The proof is similar to Construction 18. \hfill $\square$

\medskip

A \emph{weighing matrix} $W(n, k)$ of order $n$ and weight $k$ is an $n \times n$ $(0,1,-1)$-matrix, such that $WW^T = k I_n$.

\begin{construction}\label{c20}
Let $W$ be a $(4t, 4(t-1))$-weighing matrix, such that the main diagonal of $W$ contains blocks of zeros of size 4. Construct matrix $A'$ by replacing each $0$ with $O_2$, each $1$ with $I_2$ and each $-1$ with $J_2-I_2$. Then matrix $A = A' + I_t \otimes ((J_4-I_4) \otimes J_2)$ is the adjacency matrix of a DDG $\Gamma$ with parameters $(8t,4t+2,6,2t+2,4t,2)$.
\end{construction}

\proof Consider two vertices corresponding to two rows of $A$, which were obtained from one row of  $W$ after replacing. These vertices form a block of canonical partition of $\Gamma$. They have only 6 common neighbours in $\Gamma$, which come from the addition of $I_t \otimes ((J_4-I_4) \otimes J_2)$. 

Consider two vertices corresponding to two rows of $A$, which were obtained from different rows of $W$ after replacing. We have two cases: these vertices correspond to the same block of zeros from $W$ or they correspond to different blocks of zeros from $W$.

Consider two vertices from the first case. They have 4 common neighbours from addition of $I_t \otimes ((J_4-I_4) \otimes J_2)$ plus $4(t-1)/2$ common neighbours from replacing. Thus, these vertices have $2t + 2$ common neighbours.

Consider two vertices from the second case. They have 3 common neighbours in each block of $\Gamma$, which correspond to a block of zeros from $W$. Thus, they have 6 common neighbours. They also have $4(t-2)/2$ common neighbours from replacing. Thus, these vertices have $2t + 2$ common neighbours total.

So, $\Gamma$ is a DDG with parameters $(8t,4t+2,6,2t+2,4t,2)$. \hfill $\square$

\begin{construction}\label{c21}
Let $\Gamma$ be a DDG obtained from Construction 20 with adjacency matrix $A$. The main diagonal of $A$ consists of $I_t \otimes ((J_4-I_4) \otimes J_2)$, which gives a partition of $\Gamma$ into complete multipartite graphs with 4 parts of size 2. Construct a new graph $\Gamma'$ by removing the edges of the complete bipartite subgraph $K_{4,4}$ from each part of this partition. Then $\Gamma'$ is a DDG with parameters $(8t,4t-2,2,2t-2,4t,2)$.
\end{construction}

\proof The proof is similar to Construction 20. \hfill $\square$

\section{Sporadic constructions}

\begin{construction}\label{c22}
The following matrix $M$ is the adjacency matrix of a DDG with parameters $(27, 8, 4, 2, 9, 3)$:

$M = \begin{bmatrix}
O & J-I & J-I & J-I & J-I & O & O & O & O\\
J-I & O & T_1 & T_2 & T_3 & J & O & O & O\\
J-I & T_1 & O & T_3 & T_2 & O & J & O & O\\
J-I & T_2 & T_3 & O & T_1 & O & O & J & O\\
J-I & T_3 & T_2 & T_1 & O & O & O & O & J\\
O & J & O & O & O & J-I & I & I & I\\
O & O & J & O & O & I & J-I & I & I\\
O & O & O & J & O & I & I & J-I & I\\
O & O & O & O & J & I & I & I & J-I
\end{bmatrix},$

where $J = J_3, I = I_3, 
T_1 = \begin{bmatrix}
1 & 0 & 0\\
0 & 0 & 1\\
0 & 1 & 0
\end{bmatrix}, 
T_2 = \begin{bmatrix}
0 & 1 & 0\\
1 & 0 & 0\\
0 & 0 & 1
\end{bmatrix}, 
T_3 = \begin{bmatrix}
0 & 0 & 1\\
0 & 1 & 0\\
1 & 0 & 0
\end{bmatrix}.$

\end{construction}

\begin{construction}\label{c23}
The following matrix $M$ is the adjacency matrix of a DDG with parameters $(28, 6, 2, 1, 7, 4)$:

$M = \begin{bmatrix}
O &   A & B & C & O &   O &   O\\
A^T & D & O & O & E &   O &   O\\
B^T & O & D & O & O &   E &   O\\
C^T & O & O & D & O &   O &   E\\
O &   E & O & O & O &   A &   C\\
O &   O & E & O & A^T & O &   F\\
O &   O & O & E & C^T & F^T & O
\end{bmatrix}, \text{ where }$

$A = \begin{bmatrix}
1 & 1 & 0 & 0\\
1 & 1 & 0 & 0\\
0 & 0 & 1 & 1\\
0 & 0 & 1 & 1
\end{bmatrix},
B = \begin{bmatrix}
1 & 1 & 0 & 0\\
0 & 0 & 1 & 1\\
1 & 1 & 0 & 0\\
0 & 0 & 1 & 1
\end{bmatrix},
C = \begin{bmatrix}
1 & 1 & 0 & 0\\
0 & 0 & 1 & 1\\
0 & 0 & 1 & 1\\
1 & 1 & 0 & 0
\end{bmatrix},$

$D = \begin{bmatrix}
0 & 1 & 1 & 0\\
1 & 0 & 0 & 1\\
1 & 0 & 0 & 1\\
0 & 1 & 1 & 0
\end{bmatrix},
E = \begin{bmatrix}
1 & 0 & 1 & 0\\
0 & 1 & 0 & 1\\
1 & 0 & 1 & 0\\
0 & 1 & 0 & 1
\end{bmatrix},
F = \begin{bmatrix}
1 & 0 & 0 & 1\\
0 & 1 & 1 & 0\\
0 & 1 & 1 & 0\\
1 & 0 & 0 & 1
\end{bmatrix}.$

\end{construction}

\medskip

\begin{proposition}\label{p7}
There exists two DDG with parameters $(32, 10, 2, 3, 8, 4)$ for which the adjacency matrix has the given structure: $M = \begin{bmatrix}
A & B\\
B^T & C
\end{bmatrix},$ where $A$ and $C$ are the adjacency matrices of the Shrikhande graph and $B$ is the incidence matrix of a symmetric divisible design with parameters $(16, 4, 0, 1, 4, 4)$.
\end{proposition}

\section{Enumeration algorithm}
\subsection{Search for feasible parameters} 

At the first step, for a fixed number of vertices $v$, we calculate all feasible parameters ($v, k, \lambda_1, \lambda_2, m, n$) of DDGs with the approach based on M.A.~Meulenberg's work \cite{M2008}:

(1) All possible numbers of classes $m$ and sizes of classes $n$ are calculated, $m \cdot n$ must be equal to $v$.

(2) $k$ runs from 3 to $v-3$, both $\lambda_1$ and $\lambda_2$ run from max$(0,\;2k-v)$ to $k$, $\lambda_1 \not= \lambda_2$.

(3) For the remaining possibilities we check the conditions for parameters and spectrum given in Section 2.

\subsection{Enumeration of quotient matrices}

Given feasible parameters $(v, k, \lambda_1, \lambda_2, m, n)$ of a DDG, we initially construct all possible quotient matrices $R$ using the following method:

(1) We generate all possible $r_{11}$, $0 \leqslant r_{11} \leqslant n-1$.

(2) For all obtained $r_{11}$ we generate all unordered  partitions of $k - r_{11}$ as sum of $m-1$ numbers. We use unordered  partitions since the permutation of $r_{1j}$, $j \geqslant 2$ does not change the final result. 

(3) We check if the generated row satisfies the equality for $R^2 = (k^2 - \lambda_2v)I_m +\lambda_2nJ_m$. If not, we reject this row. Also in case of odd $n$ all $r_{ii}$ must be even numbers.

(4)  We repeat steps (1)--(3) for the next rows, but we generate ordered partitions instead. After we generate all possible $i \geqslant 2$ rows of the quotient matrix, we leave only non-equal up to classes renumbering options.

After we generate all non-equal quotient matrices for given parameters, in case where both $k^2 - \lambda_2 v$ and $k - \lambda_1$ are squares, we also compute all possible values of $g_1$ using the equation $$\text{trace}(A) = 0 = k + (f_1 - f_2)\sqrt{k - \lambda_1} + (g_1 - g_2)\sqrt{k^2 - \lambda_2v}.$$ 

Then we find the corresponding multiplicity for the generated quotient matrix and check if the found value is among possible values of $g_1$. If not, we reject this quotient matrix.

\subsection{Constructing adjacency matrices}

We say that two partially filled adjacency matrices are equivalent if for the graphs determined by them there is an isomorphism, which keeps partition into classes. 

For a given quotient matrix R and tuple of parameters $(v, k, \lambda_1, \lambda_2, m, n)$ we generate all possible rows of the adjacency matrix, relying on the known partition into classes. For that, we use an exhaustive search of possible rows. For each new row $t$ of the adjacency matrix ($1 \leqslant t \leqslant v$) and each corresponding entry $r_{ij}$ of the quotient matrix all possible combinations of $r_{ij}$ 1s in $n$ positions are considered.

For each $t$, all partial matrices were checked for equivalence after adding all possible $t$-th rows of the adjacency matrix and only the non-equivalent ones were considered for adding the next row.

At the last step, all obtained graphs were checked for isomorphism, and only non-isomorphic graphs were left.

We used SageMath for the search of feasible parameters, the enumeration of quotient matrices and equivalence/isomorphism check during the enumeration of adjacency matrices. For the procedure of adding a new row to the adjacency matrix, a C program was used.

\section{Enumeration results}

In Table \ref{tab:1} below we enumerate the non-trivial
proper DDGs on at most 39 vertices (trivial DDGs are DDGs that can be obtained from Constructions 1-5). The column indicated by \# gives a number of non-isomorphic DDGs with specified parameters and spectrum; $v, k, \lambda_1, \lambda_2, m, n$ are the parameters; $\theta_1^{f_1}, \theta_2^{f_2}, \theta_3^{g_1}, \theta_4^{g_2}$ are the non-principal eigenvalues with corresponding multiplicities; `WR' denotes whether all graphs with specified parameters and spectrum are walk-regular, or all graphs are not walk-regular. The column `constructions' refers to the constructions from Section~4, which can be used to obtain the specified graphs.

In column \# the exclamation mark `!' indicates that the number is the exact number of non-isomorphic graphs, a
number followed by a `+' gives the number of known non-isomorphic graphs (there can be more DDGs with the given parameters), a question mark `?' indicates that the existence of DDGs with the given parameters remains unresolved. 

For the `constructions' column, `c' is short for `construction`, `p' is short for `proposition`, `DSS' means `dual Seidel switching', `PC' means `partial complement', `L(G)' denotes the line graph of G, `J(m,n)' denotes the Johnson graph, `T(n)' denotes the triangular graph (J(2,n)). For Construction 14, we note with `.' the number of the construction (11, 12 or 13) for which the switching of edges was implemented.

{\scriptsize
\begin{longtable}{|c|c|c|c|c|c|c|c|c|c|c|c|l|}
\caption{Non-trivial proper DDGs with at most 39 vertices}
\label{tab:1}\\

\multicolumn{1}{|c|}{\#} & 
\multicolumn{1}{c|}{$v$} & 
\multicolumn{1}{c|}{$k$} & 
\multicolumn{1}{c|}{$\lambda_1$} & 
\multicolumn{1}{c|}{$\lambda_1$} & 
\multicolumn{1}{c|}{$m$} & 
\multicolumn{1}{c|}{$n$} & 
\multicolumn{1}{c|}{$\theta_1^{f_1}$} & 
\multicolumn{1}{c|}{$\theta_2^{f_2}$} & 
\multicolumn{1}{c|}{$\theta_3^{g_1}$} & 
\multicolumn{1}{c|}{$\theta_4^{g_2}$} & 
\multicolumn{1}{c|}{WR} & 
\multicolumn{1}{c|}{constructions} \\ [1pt] \hline 
& & & & & & & & & & & & \\ [-0.7em]
\endfirsthead

\multicolumn{1}{|c|}{\#} & 
\multicolumn{1}{c|}{$v$} & 
\multicolumn{1}{c|}{$k$} & 
\multicolumn{1}{c|}{$\lambda_1$} & 
\multicolumn{1}{c|}{$\lambda_1$} & 
\multicolumn{1}{c|}{$m$} & 
\multicolumn{1}{c|}{$n$} & 
\multicolumn{1}{c|}{$\theta_1^{f_1}$} & 
\multicolumn{1}{c|}{$\theta_2^{f_2}$} & 
\multicolumn{1}{c|}{$\theta_3^{g_1}$} & 
\multicolumn{1}{c|}{$\theta_4^{g_2}$} & 
\multicolumn{1}{c|}{WR} & 
\multicolumn{1}{c|}{constructions} \\ [1pt] \hline 
& & & & & & & & & & & & \\ [-0.7em]
\endhead

\hline \multicolumn{13}{c}{} \\
\endfoot

\hline \multicolumn{13}{c}{} \\
\endlastfoot

1! & 8 & 4 & 0 & 2 & 4 & 2 & $2^{1}$ & $-2^{3}$ & $0^{3}$ & -- & + & $2 \times 4$-lattice [c\ref{c8}, c\ref{c11}] \\ [2pt]
1! & 12 & 5 & 0 & 2 & 6 & 2 & $\sqrt{5}^{3}$ & $-\sqrt{5}^{3}$ & -- & $-1^{5}$ & + & icosahedron [c\ref{c6}] \\ [2pt]
1! & 12 & 5 & 1 & 2 & 4 & 3 & $2^{2}$ & $-2^{6}$ & $1^{3}$ & -- & + & $3 \times 4$-lattice [c\ref{c8}, c\ref{c11}] \\ [2pt]
1! & 12 & 6 & 2 & 3 & 3 & 4 & $2^{3}$ & $-2^{6}$ & $0^{2}$ & -- & + & L(octahedron) \\ [2pt]
1! & 12 & 7 & 3 & 4 & 4 & 3 & $2^{2}$ & $-2^{6}$ & $1^{2}$ & $-1^{1}$ & + & c\ref{c8}, c\ref{c14}.\ref{c11} \\ [2pt]
1! & 15 & 4 & 0 & 1 & 5 & 3 & $2^{5}$ & $-2^{5}$ & -- & $-1^{4}$ & + & L(Petersen graph) [c\ref{c6}] \\ [2pt]
1! & 18 & 9 & 6 & 4 & 6 & 3 & $\sqrt{3}^{6}$ & $-\sqrt{3}^{6}$ & $3^{1}$ & $-3^{4}$ & + & c\ref{c18} from $(9, 3, 0, 1)$-SDD \\ [2pt]
1! & 20 & 7 & 3 & 2 & 4 & 5 & $2^{4}$ & $-2^{12}$ & $3^{3}$ & -- & + & $5 \times 4$-lattice [c\ref{c8}, c\ref{c11}] \\ [2pt]
1! & 20 & 9 & 0 & 4 & 10 & 2 & $3^{4}$ & $-3^{6}$ & $1^{3}$ & $-1^{6}$ & -- & DSS(J(6,3)) \\ [2pt]
1! & 20 & 9 & 0 & 4 & 10 & 2 & $3^{5}$ & $-3^{5}$ & -- & $-1^{9}$ & + & J(6,3) [c\ref{c6}] \\ [2pt]
1! & 20 & 13 & 9 & 8 & 4 & 5 & $2^{4}$ & $-2^{12}$ & $3^{2}$ & $-3^{1}$ & + & c\ref{c8}, c\ref{c14}.\ref{c11} \\ [2pt]
1! & 24 & 6 & 2 & 1 & 3 & 8 & $2^{9}$ & $-2^{12}$ & $\sqrt{12}^{1}$ & $-\sqrt{12}^{1}$ & -- & c\ref{c9} \\ [2pt]
1! & 24 & 7 & 0 & 2 & 8 & 3 & $\sqrt{7}^{8}$ & $-\sqrt{7}^{8}$ & -- & $-1^{7}$ & + & Klein graph [c\ref{c6}] \\ [2pt]
1! & 24 & 8 & 4 & 2 & 4 & 6 & $2^{5}$ & $-2^{15}$ & $4^{3}$ & -- & + & $6 \times 4$-lattice [c\ref{c8}, c\ref{c11}] \\ [2pt]
1! & 24 & 8 & 4 & 2 & 4 & 6 & $2^{7}$ & $-2^{13}$ & $4^{2}$ & $-4^{1}$ & -- & c\ref{c12} \\ [2pt]
6! & 24 & 8 & 4 & 2 & 4 & 6 & $2^{9}$ & $-2^{11}$ & $4^{1}$ & $-4^{2}$ & + & c\ref{c13} \\ [2pt]
5! & 24 & 10 & 2 & 4 & 12 & 2 & $\sqrt{8}^{6}$ & $-\sqrt{8}^{6}$ & $2^{3}$ & $-2^{8}$ & + & c\ref{c21} \\ [2pt]
2! & 24 & 10 & 3 & 4 & 8 & 3 & $\sqrt{7}^{8}$ & $-\sqrt{7}^{8}$ & $2^{1}$ & $-2^{6}$ & + & one graph is Cayley (see \cite{GS2014}) \\ [2pt]
1! & 24 & 10 & 6 & 3 & 3 & 8 & $2^{8}$ & $-2^{13}$ & $\sqrt{28}^{1}$ & $-\sqrt{28}^{1}$ & -- & c\ref{c9} \\ [2pt]
1! & 24 & 14 & 6 & 8 & 12 & 2 & $\sqrt{8}^{6}$ & $-\sqrt{8}^{6}$ & $2^{2}$ & $-2^{9}$ & + & c\ref{c10}, c\ref{c20} \\ [2pt]
1! & 24 & 14 & 7 & 8 & 8 & 3 & $\sqrt{7}^{8}$ & $-\sqrt{7}^{8}$ & -- & $-2^{7}$ & + & PC(Klein graph) [p\ref{p5}] \\ [2pt]
1! & 24 & 16 & 12 & 10 & 4 & 6 & $2^{5}$ & $-2^{15}$ & $4^{2}$ & $-4^{1}$ & + & c\ref{c8}, c\ref{c14}.\ref{c11} \\ [2pt]
1! & 24 & 16 & 12 & 10 & 4 & 6 & $2^{7}$ & $-2^{13}$ & $4^{1}$ & $-4^{2}$ & -- & c\ref{c14}.\ref{c12} \\ [2pt]
4! & 24 & 16 & 12 & 10 & 4 & 6 & $2^{9}$ & $-2^{11}$ & -- & $-4^{3}$ & + & c\ref{c14}.\ref{c13} \\ [2pt]
1! & 27 & 8 & 4 & 2 & 9 & 3 & $2^{7}$ & $-2^{11}$ & $\sqrt{10}^{4}$ & $-\sqrt{10}^{4}$ & -- & c\ref{c22} \\ [2pt]
2! & 27 & 18 & 9 & 12 & 9 & 3 & $3^{6}$ & $-3^{12}$ & $0^{8}$ & -- & + & c\ref{c16} from Schl\"{a}fli graph \\ [2pt]
1! & 28 & 6 & 2 & 1 & 7 & 4 & $2^{9}$ & $-2^{12}$ & $\sqrt{8}^{3}$ & $-\sqrt{8}^{3}$ & -- & c\ref{c23} \\ [2pt]
1! & 28 & 9 & 5 & 2 & 4 & 7 & $2^{6}$ & $-2^{18}$ & $5^{3}$ & -- & + & $7 \times 4$-lattice [c\ref{c8}, c\ref{c11}] \\ [2pt]
1! & 28 & 13 & 0 & 6 & 14 & 2 & $\sqrt{13}^{7}$ & $-\sqrt{13}^{7}$ & -- & $-1^{13}$ & + & Taylor graph [c\ref{c6}] \\ [2pt]
16! & 28 & 13 & 4 & 6 & 7 & 4 & $3^{9}$ & $-3^{12}$ & $1^{1}$ & $-1^{5}$ & -- & -- \\ [2pt]
56! & 28 & 15 & 6 & 8 & 7 & 4 & $3^{7}$ & $-3^{14}$ & $1^{6}$ & -- & + & c\ref{c16} from T(8), Chang graphs \\ [2pt]
4! & 28 & 15 & 6 & 8 & 7 & 4 & $3^{8}$ & $-3^{13}$ & $1^{3}$ & $-1^{3}$ & -- & DSS of the previous entry \\ [2pt]
1! & 28 & 19 & 15 & 12 & 4 & 7 & $2^{6}$ & $-2^{18}$ & $5^{2}$ & $-5^{1}$ & + & c\ref{c8}, c\ref{c14}.\ref{c11} \\ [2pt]
2! & 32 & 10 & 2 & 3 & 8 & 4 & $\sqrt{8}^{12}$ & $-\sqrt{8}^{12}$ & $2^{1}$ & $-2^{6}$ & + & p\ref{p7} \\ [2pt]
1! & 32 & 10 & 6 & 2 & 4 & 8 & $2^{7}$ & $-2^{21}$ & $6^{3}$ & -- & + & $8 \times 4$-lattice [c\ref{c8}, c\ref{c11}] \\ [2pt]
1! & 32 & 10 & 6 & 2 & 4 & 8 & $2^{10}$ & $-2^{18}$ & $6^{2}$ & $-6^{1}$ & -- & c\ref{c12} \\ [2pt]
15! & 32 & 10 & 6 & 2 & 4 & 8 & $2^{13}$ & $-2^{15}$ & $6^{1}$ & $-6^{2}$ & + & c\ref{c13} \\ [2pt]
15! & 32 & 14 & 2 & 6 & 16 & 2 & $\sqrt{12}^{8}$ & $-\sqrt{12}^{8}$ & $2^{4}$ & $-2^{11}$ & + & c\ref{c21} \\ [2pt]
2+ & 32 & 15 & 6 & 7 & 4 & 8 & $3^{12}$ & $-3^{16}$ & -- & $-1^{3}$ & + & Cayley graphs (see \cite{GS2014}) \\ [2pt]
1! & 32 & 16 & 0 & 8 & 16 & 2 & $4^{6}$ & $-4^{10}$ & $0^{15}$ & -- & + & c\ref{c15} from halved 6-cube \\ [2pt]
? & 32 & 17 & 8 & 9 & 4 & 8 & $3^{17}$ & $-3^{11}$ & $1^{2}$ & $-1^{1}$ & -- & -- \\ [2pt]
1! & 32 & 18 & 6 & 10 & 16 & 2 & $\sqrt{12}^{8}$ & $-\sqrt{12}^{8}$ & $2^{3}$ & $-2^{12}$ & + & c\ref{c20} \\ [2pt]
1! & 32 & 22 & 18 & 14 & 4 & 8 & $2^{7}$ & $-2^{21}$ & $6^{2}$ & $-6^{1}$ & + & c\ref{c8}, c\ref{c14}.\ref{c11} \\ [2pt]
1! & 32 & 22 & 18 & 14 & 4 & 8 & $2^{10}$ & $-2^{18}$ & $6^{1}$ & $-6^{2}$ & -- & c\ref{c14}.\ref{c12} \\ [2pt]
9! & 32 & 22 & 18 & 14 & 4 & 8 & $2^{13}$ & $-2^{15}$ & -- & $-6^{3}$ & + & c\ref{c14}.\ref{c13} \\ [2pt]
2! & 35 & 12 & 3 & 4 & 7 & 5 & $3^{12}$ & $-3^{16}$ & $2^{3}$ & $-2^{3}$ & -- & -- \\ [2pt]
3854! & 35 & 12 & 3 & 4 & 7 & 5 & $3^{14}$ & $-3^{14}$ & -- & $-2^{6}$ & + & c\ref{c7} from $(35, 18, 9)$-graphs \\ [2pt]
3! & 36 & 9 & 3 & 2 & 12 & 3 & $\sqrt{6}^{12}$ & $-\sqrt{6}^{12}$ & $3^{4}$ & $-3^{7}$ & + & two graphs from c\ref{c19} \\ [2pt]
7! & 36 & 9 & 4 & 2 & 18 & 2 & $\sqrt{5}^{9}$ & $-\sqrt{5}^{9}$ & $3^{7}$ & $-3^{10}$ & + & c\ref{c17} from icosahedron \\ [2pt]
1! & 36 & 11 & 7 & 2 & 4 & 9 & $2^{8}$ & $-2^{24}$ & $7^{3}$ & -- & + & $8 \times 4$-lattice [c\ref{c8}, c\ref{c11}] \\ [2pt]
1! & 36 & 17 & 0 & 8 & 18 & 2 & $\sqrt{17}^{9}$ & $-\sqrt{17}^{9}$ & -- & $-1^{17}$ & + & Taylor graph [c\ref{c6}] \\ [2pt]
3+ & 36 & 24 & 15 & 16 & 4 & 9 & $3^{12}$ & $-3^{20}$ & $0^{3}$ & -- & + & Cayley graphs (see \cite{GS2014}) \\ [2pt]
1! & 36 & 25 & 21 & 16 & 4 & 9 & $2^{8}$ & $-2^{24}$ & $7^{2}$ & $-7^{1}$ & + & c\ref{c8}, c\ref{c14}.\ref{c11} \\ [2pt]
1! & 36 & 27 & 21 & 20 & 12 & 3 & $\sqrt{6}^{12}$ & $-\sqrt{6}^{12}$ & $3^{1}$ & $-3^{10}$ & + & c\ref{c18} from $(18, 12, 6, 8)$-SDD \\ [2pt]
2! & 38 & 9 & 0 & 2 & 19 & 2 & $3^{8}$ & $-3^{11}$ & $\sqrt{5}^{9}$ & $-\sqrt{5}^{9}$ & -- & -- \\ [2pt]

\end{longtable}}

The list of proper divisible design graphs up to 39 vertices is available by \url{http://alg.imm.uran.ru/dezagraphs/ddgtab.html}. This web page provides access to adjacency matrices, quotient matrices and other properties of the graphs we found.

\section{Conclusion}

The enumeration was incomplete for three tuples of parameters: $(32,15,6,7,4,8)$, $(32,17,8,9,4,8)$, $(36,24,15,16,4,9)$. For these sets, the adjacency matrix construction algorithm was enumerating the rows of adjacency matrices corresponding to the first row of the quotient matrix for several weeks. Therefore, we decided to stop the enumeration for these tuples and switch to others. It is possible to enumerate graphs for these tuples with more time or more powerful equipment and it could be done in the future.

After the enumeration was finished, we compared our results with the results from \cite{CH2011, CH2014, HKM2011}. We found that the results for the existence of graphs coincide, except for one tuple of parameters: $(27, 8, 4, 2, 9, 3)$. In \cite{HKM2011} this tuple was rejected because it did not meet the necessary condition given in \cite[Theorem 5.1]{HKM2011}, which concerns the existence of non-zero integral solutions of the Diophantine equations. It turns out that the Diophantine equations for parameters $(27, 8, 4, 2, 9, 3)$ actually have required solutions, and our enumeration produced one graph with these parameters. Our enumeration also showed the non-existence of graphs with parameters $(27, 16, 12, 9, 9, 3)$, the only case for which the answer was not given in \cite{CH2011, CH2014, HKM2011}.

For one graph with parameters $(24, 10, 3, 4, 8, 3)$, sixteen graphs with parameters $(28, 13, 4, 6, 7, 4)$, two graphs with parameters $(35, 12, 3, 4, 7, 5)$, one graph with parameters $(36, 9, 3, 2, 12, 3)$ and two graphs with parameters $(38, 9, 0, 2, 19, 2)$ we could not find existing theoretical constructions or present new ones. We hope that one day new theoretical constructions will be found for the graphs that remained undescribed.

\section*{Acknowledgments}
The authors are grateful to the referees for their valuable remarks and suggestions which improved the paper.

The reported study was funded by RFBR according to the research project 20-51-53023.

\end{document}